\documentclass[
  a4paper,oneside,DIV=12,
  12pt,
  headsepline,
  ]{scrartcl}

\usepackage{etoolbox}
\usepackage{ifdraft}
\usepackage{iftex}
\usepackage{subfiles}

\ifluatex
  \usepackage[utf8]{luainputenc}
\else
  \usepackage[utf8]{inputenc}
  \usepackage[T1]{fontenc}
\fi

\usepackage[english]{babel}
\usepackage[autostyle=true]{csquotes}

\usepackage{lmodern}
\usepackage{fourier}

\usepackage{amsfonts}
\usepackage{mathrsfs} %
\usepackage{dsfont}

\usepackage{amssymb}

\usepackage{amsmath}
\usepackage{mathtools}
\usepackage[thmmarks, amsmath, amsthm]{ntheorem}
\usepackage{nccmath}    %
\usepackage{tensor}

\usepackage[
  hscale=0.8,vscale=0.875,
  includehead,footskip=2\baselineskip,
]{geometry}

\usepackage[draft=false]{scrlayer-scrpage}
\usepackage{needspace}

\usepackage[cmyk,dvipsnames]{xcolor}

\ifdraft{%
  \usepackage[textsize=scriptsize,colorinlistoftodos]{todonotes}
  \usepackage{lineno}
  \usepackage{scrtime}
}{}

\usepackage{booktabs}
\usepackage[inline]{enumitem}
\usepackage{lettrine}
\usepackage{url}

\usepackage[backend=biber,style=numeric-comp,giveninits,url=false,sortcites=true,maxbibnames=99]{biblatex}
\usepackage[final,
  pdfborder={0 0 0}, colorlinks=true,
  linkcolor=BrickRed, citecolor=ForestGreen, urlcolor=RoyalBlue]{hyperref}
\usepackage{zref-clever}

\ifdraft{%
  \linenumbers
  \newcommand*\patchAmsMathEnvironmentForLineno[1]{%
    \expandafter\let\csname old#1\expandafter\endcsname\csname #1\endcsname
    \expandafter\let\csname oldend#1\expandafter\endcsname\csname end#1\endcsname
    \renewenvironment{#1}%
    {\linenomath\csname old#1\endcsname}%
    {\csname oldend#1\endcsname\endlinenomath}}%
  \newcommand*\patchBothAmsMathEnvironmentsForLineno[1]{%
    \patchAmsMathEnvironmentForLineno{#1}%
    \patchAmsMathEnvironmentForLineno{#1*}}%
  \AtBeginDocument{%
    \patchBothAmsMathEnvironmentsForLineno{equation}%
    \patchBothAmsMathEnvironmentsForLineno{align}%
    \patchBothAmsMathEnvironmentsForLineno{flalign}%
    \patchBothAmsMathEnvironmentsForLineno{alignat}%
    \patchBothAmsMathEnvironmentsForLineno{gather}%
    \patchBothAmsMathEnvironmentsForLineno{multline}%
  }}

\clearmainofpairofpagestyles
\automark[section]{subsection}

\renewcommand{\subsectionmark}[1]{}

\lohead{{\small\normalfont \headertitle}}
\rohead{{\small \headerauthors}}

\RedeclareSectionCommand[%
  font=\Large\sffamily\bfseries,%
  beforeskip=1\baselineskip,%
  afterskip=0.5\baselineskip,%
  indent=0em,%
  afterindent=false,%
  tocbeforeskip=0.3\baselineskip plus 1pt minus 1pt%
]{section}

\RedeclareSectionCommands[%
  font=\normalfont\bfseries,%
  beforeskip=3pt,%
  afterskip=-1em%
]{subsection,subsubsection}

\RedeclareSectionCommands[%
  font=\normalfont\itshape,%
  beforeskip=2pt,%
  afterskip=-1em,%
  indent=0pt,%
]{paragraph}

\setitemize{itemsep=0.02\baselineskip}

\newenvironment{enumeratearabic}{
  \begin{enumerate}[label=(\arabic*),%
  leftmargin=2.5em,itemindent=0pt,%
  labelindent=.5em,labelwidth=1.5em,labelsep=!,%
  noitemsep]
  }{
  \end{enumerate}
}

\newenvironment{enumeratearabic*}{
  \begin{enumerate*}[label=(\arabic*)] %
    }{
  \end{enumerate*}
}

\newenvironment{enumerateroman*}{
  \begin{enumerate*}[label=(\roman*)] %
    }{
  \end{enumerate*}
}

\numberwithin{equation}{section}

\theoremnumbering{arabic}
\newtheorem{theoremcounter}{theoremcounter}[section]
\theoremnumbering{Alph}
\newtheorem{maintheoremcounter}{maintheoremcounter}

\theoremstyle{plain}

\newtheorem{lemma}[theoremcounter]{Lemma}

\newtheorem{proposition}[theoremcounter]{Proposition}

\theoremstyle{plain}

\newtheorem{maintheorem}[maintheoremcounter]{Theorem}

\theoremstyle{definition}

\theoremstyle{remark}

\newtheorem{remark}[theoremcounter]{Remark}

\theoremstyle{nonumberremark}

\newtheorem{mainremark}{Remark}

\newenvironment{mainremarkenumerate}
{%
  \mainremark
  \enumeratearabic
}{%
  \endenumeratearabic
  \endmainremark
}%

\zcsetup{
  cap = true,
}

\zcLanguageSetup{english}{
  type = maintheorem ,
    Name-sg = Theorem ,  name-sg = theorem ,
    Name-pl = Theorems , name-pl = theorems ,
  type = maincorollary ,
    Name-sg = Corollary ,  name-sg = corollary ,
    Name-pl = Corollaries , name-pl = corollaries ,
  type = mainremark ,
    Name-sg = Remark ,  name-sg = remark ,
    Name-pl = Remarks , name-pl = remarks ,
  type = assumption ,
    Name-sg = Assumption ,  name-sg = assumption ,
    Name-pl = Assumptions , name-pl = assumptions ,
  type = claim ,
    Name-sg = Claim ,  name-sg = claim ,
    Name-pl = Claims , name-pl = claims ,
  type = conjecture ,
    Name-sg = Conjecture ,  name-sg = conjecture ,
    Name-pl = Conjectures , name-pl = conjectures ,
  type = notation ,
    Name-sg = Notation ,  name-sg = notation ,
    Name-pl = Notations , name-pl = notations ,
  type = problem ,
    Name-sg = Problem ,  name-sg = problem ,
    Name-pl = Problems , name-pl = problems ,
  type = question ,
    Name-sg = Question ,  name-sg = question ,
    Name-pl = Questions , name-pl = questions ,
}

\newcommand{\zctheoremtype}[2]{%
  \AddToHook{env/#1/begin}{\zcsetup{countertype = {#2 = #1}}}%
}

\zctheoremtype{theorem}{theoremcounter}
\zctheoremtype{lemma}{theoremcounter}
\zctheoremtype{corollary}{theoremcounter}
\zctheoremtype{proposition}{theoremcounter}
\zctheoremtype{definition}{theoremcounter}
\zctheoremtype{remark}{theoremcounter}
\zctheoremtype{example}{theoremcounter}
\zctheoremtype{exercise}{theoremcounter}
\zctheoremtype{assumption}{theoremcounter}
\zctheoremtype{claim}{theoremcounter}
\zctheoremtype{conjecture}{theoremcounter}
\zctheoremtype{notation}{theoremcounter}
\zctheoremtype{problem}{theoremcounter}
\zctheoremtype{question}{theoremcounter}
\zctheoremtype{maintheorem}{maintheoremcounter}
\zctheoremtype{maincorollary}{maintheoremcounter}

\AtEveryBibitem{\clearfield{doi}}
\AtEveryBibitem{\clearfield{isbn}}
\AtEveryBibitem{\clearfield{issn}}
\AtEveryBibitem{\clearfield{pages}}
\AtEveryBibitem{\clearlist{language}}

\renewcommand*{\bibfont}{\footnotesize}
\renewbibmacro*{in:}{}
\DeclareFieldFormat
[article,inbook,incollection,inproceedings,patent,thesis,unpublished,misc]
{title}{#1}

\newcommand{\tx}{\text}
\newcommand{\nbd}{\nobreakdash-\hspace{0pt}}
\newcommand{\thdash}{\nbd th}

\newcommand{\fref}[2]{\hyperref[#2]{#1~\ref*{#2}}}

\makeatletter
\newcommand{\writelabel}[1]{#1\def\@currentlabel{#1}}
\makeatother

\newcommand{\minwidthmathbox}[2]{%
  \mathmakebox[{\ifdim#1<\width\width\else#1\fi}]{#2}%
}

\newcommand{\tbf}{\bfseries}

\newcommand{\td}{\tilde}
\newcommand{\wtd}{\widetilde}
\newcommand{\ov}{\overline}

\newcommand{\bbM}{\mathbb{M}}

\newcommand{\frakd}{\mathfrak{d}}
\newcommand{\frake}{\mathfrak{e}}

\newcommand{\rmM}{\mathrm{M}}

\newcommand{\rmS}{\mathrm{S}}

\newcommand{\defcol}{\mathrel{:}}
\newcommand{\defeq}{\mathrel{:=}}

\newcommand{\condsep}{\mathrel{\;:\;}}

\newcommand{\mrelspace}[1]{\mathrel{\mspace{#1}}}

\NewCommandCopy{\rightarroworig}{\rightarrow}
\renewcommand{\rightarrow}{%
  \protect\relbar\mrelspace{-9.7mu}\rightarroworig}

\NewCommandCopy{\leftarroworig}{\leftarrow}
\renewcommand{\leftarrow}{%
  \protect\leftarroworig\mrelspace{-9.7mu}\relbar}

\newcommand{\ra}{\rightarrow}

\NewCommandCopy{\slashorig}{\slash}
\renewcommand{\slash}{\mathbin{\slashorig}}

\renewcommand{\pmod}[1]{\;(\mathrm{mod}\, #1)}

\newcommand{\sgn}{\mathrm{sgn}}

\newenvironment{smatrix}{\begin{smallmatrix}}{\end{smallmatrix}}
\newenvironment{psmatrix}{\left(\begin{smatrix}}{\end{smatrix}\right)}

\newcommand{\ZZ}{\mathbb{Z}}

\newcommand{\CC}{\mathbb{C}}

\newcommand{\Mp}[1]{\mathrm{Mp}_{#1}}

\newcommand{\SL}[1]{\mathrm{SL}_{#1}}

\newcommand{\U}[1]{\mathrm{U}_{#1}}

\newcommand{\HS}{\mathbb{H}}

\newcommand{\Ga}{\Gamma}

\newcommand{\hol}{\mathrm{hol}}

\newcommand{\rhoThree}{\rho_3}

\title{%
  Recursions for Mock Theta Functions
}
\subtitle{%
  Enabling AI-Assisted Research\\for the Broader Mathematical Community
}
\author{%
  Matthew Ortiz%
  \and%
  Martin Raum%
  \thanks{The author was partially supported by Vetenskapsr\aa det Grant~2023-04217.}%
  \and%
  Olav K. Richter
  \thanks{The author was partially supported by a grant from the Simons Foundation (\#835652 to Olav Richter).}
}
\newcommand{\headertitle}{%
  Weighted Recursions for Mock Theta Functions
}
\newcommand{\headerauthors}{%
  M.~Ortiz,
  M.~Raum,
  O.~Richter
}
\ifdraft{\date{Draft: \today\ at\ \thistime}}{\date{}}

\begin{document}

\thispagestyle{scrplain}
\begingroup
\deffootnote[1em]{1.5em}{1em}{\thefootnotemark}
\maketitle
\endgroup

\begin{abstract}
\small
\noindent
{\tbf Abstract:}
We establish weighted recursions for the coefficients of Ramanujan's third order mock theta functions~$f$ and~$\omega$.
Specifically, we apply a holomorphic projection operator to vector-valued Rankin–Cohen brackets of completed mock theta series and their shadows.
By employing a vector-valued framework, we exploit the vanishing of certain spaces of vector-valued cusp forms.
Our proof is AI-assisted and prioritizes accessibility, allowing for straightforward customization and replication within the broader research community.
\\[.3\baselineskip]
\textsf{\textbf{%
    mock theta functions%
  }}%
~{\tiny$\blacksquare$}\ %
\textsf{\textbf{%
    vector-valued mock modular forms%
  }}%
~{\tiny$\blacksquare$}\ %
\textsf{\textbf{%
    holomorphic projection%
  }}%
~{\tiny$\blacksquare$}\ %
\textsf{\textbf{%
    Rankin--Cohen brackets%
  }}
\\[.2\baselineskip]
\noindent
\textsf{\textbf{%
    MSC Primary: 11F37%
  }}%
~{\tiny$\blacksquare$}\ %
\textsf{\textbf{%
    MSC Secondary: 11F27, 11F30, 68T01%
  }}
\end{abstract}

\Needspace*{4em}

\Needspace*{4em}
\phantomsection
\label{sec:introduction}
\addcontentsline{toc}{section}{Introduction}
\markright{Introduction}

\lettrine[lines=2,nindent=.2em]{\tbf M}{\,ock} theta functions have a rich history with deep connections to different areas of mathematics and physics.
For overviews of the subject, see \cite{ono-2009, zagier-2009}.
Of particular interest have been the Fourier coefficients of mock theta functions.
In~\cite{imamoglu-raum-richter-2014}, the authors determined finite recursions for the Fourier coefficients of the third order mock theta functions
\begin{gather*}
  f(q)
  \defeq
  \sum_{n=0}^{\infty} c_f(n)\, q^n
  \defeq
  1 + \sum_{n=1}^{\infty} \frac{q^{n^2}}{(1 + q)^2\, (1 + q^2)^2\,  \cdots\,  (1 + q^n)^2}
  \tx{,}
\end{gather*}
and
\begin{gather*}
  \omega(q)
  \defeq
  \sum_{n=0}^{\infty} c_{\omega}(n)\, q^n
  \defeq
  \sum_{n=0}^{\infty} \frac{q^{2n^2 + 2n}}{(1 - q)^2\, (1 - q^3)^2\,  \cdots\,  (1 - q^{2n+1})^2}
  \tx{.}
\end{gather*}
The proofs of the relations in~\cite{imamoglu-raum-richter-2014} relied on applying a holomorphic projection operator to mock modular forms.
More recently, in~\cite{ortiz-raum-richter-2026-preprint} we established new recursions for Hurwitz class numbers with polynomial weights.
Our proofs used Rankin-Cohen brackets together with a holomorphic projection operator to mock modular forms.
The difference between~\cite{ortiz-raum-richter-2026-preprint} and prior related work was the genuinely vector-valued approach, which allowed us to leverage the vanishing of spaces of vector-valued cusp forms, in contrast to the scalar case.

In this paper, we combine ideas from~\cite{imamoglu-raum-richter-2014} and~\cite{ortiz-raum-richter-2026-preprint} to provide an AI-assisted proof of new weighted recursions for the mock theta functions $f$ and $\omega$.
Specifically, we establish the following theorem, where we use the conventions that $c_f(n)=c_{\omega}(n)=0$ if $n\not\in\ZZ_{\ge0}$.
For~$s\in\ZZ_{\ge0}$, we set~$\sigma_s(x)=\sum_{d\mid x}d^s$ if~$x\in\ZZ_{>0}$ and~$\sigma_s(x)=0$ otherwise.

\begin{maintheorem}%
\label{mainthm:mock_theta_recursions}
For~$N \in \frac{1}{2}\ZZ$ and~$m\in\ZZ$ set
\begin{gather*}
  g_3(N,m)
  \defeq
  m^3-18mN
  \tx{,}\quad
  g_5(N,m)
  \defeq
  m^5-30m^3N+180mN^2
  \tx{,}
\end{gather*}
and for~$N \in \frac{1}{2}\ZZ$ and~$\nu \in \{1,2\}$ set
\begin{gather*}
  \lambda^f_\nu(N)
  \defeq
  \mspace{-30mu}
  \sum_{\substack{
      d,e\in\ZZ\\
      de=24N\\
      \frac{e-d}{2},\,\frac{e+d}{2}\in6\ZZ+1
    }}
  \mspace{-24mu}
  \sgn(e^2-d^2)\,
  \min\bigl\{ |d|,|e| \bigr\}^{\nu}
  \tx{,}
  \quad
  \lambda^\omega_\nu(N)
  \defeq
  \mspace{-30mu}
  \sum_{\substack{
      d,e\in\ZZ\\
      de=24N\\
      \frac{e-d}{2},\,\frac{e+d}{2}\in6\ZZ+2
    }}
  \mspace{-24mu}
  (-1)^{\frac{e-4}{6}}\,
  \sgn(e^2-d^2)\,
  \min\bigl\{ |d|,|e| \bigr\}^{\nu}
  \tx{.}
\end{gather*}
Then for all~$N\in\ZZ_{>0}$ the coefficients of the mock theta function~$f$ satisfy
\begin{alignat*}{3}
  2
  &
  \sum_{m\in6\ZZ+1}
  \mspace{-8mu}
  &&
  g_3(N,m)\,
  c_f\bigl(N-\mfrac{m^2-1}{24}\bigr)
  && =
  \lambda^f_3(N)
  - 32
  \bigl(
  \sigma_3(N)
  -
  16\, \sigma_3\bigl( \tfrac{N}{2} \bigr)
  \bigr)
  \tx{,}
  \\
  8
  &
  \sum_{m\in6\ZZ+1}
  \mspace{-8mu}
  &&
  g_5(N,m)\,
  c_f\bigl(N-\mfrac{m^2-1}{24}\bigr)
  && =
  \lambda^f_5(N)
  +
  64 \bigl(
  \sigma_5(N)
  -
  64\, \sigma_5\bigl( \tfrac{N}{2} \bigr)
  \bigr)
  \tx{,}
\intertext{
and for~$N\in\frac{1}{2}\ZZ_{>0}$ those of\/~$\omega$ satisfy
}
  4
  &
  \sum_{m\in6\ZZ+2}
  \mspace{-8mu}
  (-1)^{\frac{m-2}{6}}\,
  &&
  g_3(N,m)\,
  c_{\omega}\bigl(2N-\mfrac{m^2+8}{12}\bigr)
  && =
  -
  \lambda^\omega_3(N)
  -
  32
  \bigl(
  \sigma_3(2N)
  -
  \sigma_3(N)
  \bigr)
  \tx{,}
  \\
  16
  &
  \sum_{m\in6\ZZ+2}
  \mspace{-8mu}
  (-1)^{\frac{m-2}{6}}\,
  &&
  g_5(N,m)\,
  c_{\omega}\bigl(2N-\mfrac{m^2+8}{12}\bigr)
  && =
  -
  \lambda^\omega_5(N)
  +
  64
  \bigl(
  \sigma_5(2N)
  -
  \sigma_5(N)
  \bigr)
  \tx{.}
\end{alignat*}
\end{maintheorem}

\begin{mainremarkenumerate}
\item
The functions that yield the results in \zcref{mainthm:mock_theta_recursions} belong to a~$9$\nbd{}dimensional $\CC$\nbd{}vector space.
The first pair of relations arises from the first component, while the second pair is derived from the fifth component; see \zcref{ssec:holomorphic_projection_fourier_coefficients:diagonal_coefficients} and \zcref{sec:proof_main_theorem}.
\item
The machinery of this paper, along with the AI-assisted methodology, can also be applied to derive new recursions for generalized Hurwitz class numbers, as well as Fourier coefficients of other mock theta functions.
\end{mainremarkenumerate}

The proof of \zcref{mainthm:mock_theta_recursions} combines the tools from~\cite{imamoglu-raum-richter-2014} and~\cite{ortiz-raum-richter-2026-preprint}: We recall Zwegers's~\cite{zwegers-thesis, zwegers-2001} completion of mock theta functions into real-analytic vector-valued modular forms, and then apply holomorphic projection to the Rankin-Cohen brackets of Zwegers's completions and their shadows (theta functions).
As in~\cite{ortiz-raum-richter-2026-preprint}, we take advantage of the vanishing of certain spaces of vector-valued cusp forms (in weights~$4$ and~$6$) to prove our results.

The details were worked out by a \emph{large language model}, popularly known as \emph{AI}.
Specifically, we employed the model GPT\nbd{}5.5 at \enquote{xhigh} reasoning effort via a standard subscription, using merely~12 prompts.
Our approach distinguishes itself from a naive use of such models by dedicated \emph{agent skills}, which we make available.
We document the AI assistance in detail in \zcref{sec:ai-contribution}.
The mathematical methodology there is fully explicit rather than hidden inside a complex AI model, providing readers with verifiable instructions that can be customized and replicated.

The human authors of this paper verified all proofs.

\section{Preliminaries}%
\label{sec:preliminaries}
\subsection{The metaplectic group and vector-valued forms}%
\label{ssec:preliminaries:metaplectic}

We introduce some standard notation.
Let~$\HS \defeq \{ \tau = x + i y \in \CC \condsep y > 0\}$ be the Poincar\'e upper half plane.
For~$z \in \CC$, set~$e(z) \defeq \exp(2 \pi i\, z)$, for~$\tau\in\HS$, let~$q \defeq e(\tau)$, and for an integer~$r$ put~$\zeta_r \defeq e(1 / r)$.
The metaplectic group~$\Mp{2}(\ZZ)$ is generated by
\begin{gather*}
  S
  =
  \bigl(
  \begin{psmatrix} 0 & -1 \\ 1 & 0 \end{psmatrix},
  \sqrt{\tau}
  \bigr)
  \quad\tx{and}\quad
  T
  =
  \bigl(
  \begin{psmatrix} 1 & 1 \\ 0 & 1 \end{psmatrix},
  1
  \bigr)
  \tx{,}
\end{gather*}
where the square root denotes the principal branch satisfying~$\sqrt{i} =\zeta_8$.
For a finite-dimensional unitary representation~$\rho$ of~$\Mp{2}(\ZZ)$, the space~$\rmM_k(\rho)$ consists of holomorphic functions~$F\defcol\HS\ra V(\rho)$ that are invariant under the weight~$k$ slash action of type~$\rho$ and are bounded at the cusp.
We write~$\rmS_k(\rho)$ for the space of cusp forms in~$\rmM_k(\rho)$ and we denote the space of harmonic weak Maass forms of weight~$k$ and type~$\rho$ by~$\bbM_k(\rho)$; for more details see also~\cite{bringmann-folsom-ono-rolen-2018, imamoglu-raum-richter-2014,ortiz-raum-richter-2026-preprint}.

\subsection{Zwegers's vector-valued completions of $f$ and $\omega$}%
\label{ssec:preliminaries:mock_theta_vector}

Set
\begin{gather*}
  F^+(\tau)
  \defeq
  \begin{pmatrix}
    q^{-\frac{1}{24}} f(q)                    \\
    2 q^{\frac{1}{3}} \omega(q^{\frac{1}{2}}) \\
    2 q^{\frac{1}{3}} \omega(-q^{\frac{1}{2}})
  \end{pmatrix}
  \tx{,}
\end{gather*}
and~$F^-(\tau) \defeq -2 \sqrt{6} G^*(\tau)$ with
\begin{align}
  \label{eq:definition-of-G}
  G(\tau)
  \defeq
  \begin{pmatrix}
    -\sum_{n \in \ZZ} (n+\frac{1}{6}) q^{\frac{3}{2}(n+\frac{1}{6})^2}       \\
    \sum_{n \in \ZZ} (-1)^n (n+\frac{1}{3}) q^{\frac{3}{2}(n+\frac{1}{3})^2} \\
    -\sum_{n \in \ZZ} (n+\frac{1}{3}) q^{\frac{3}{2}(n+\frac{1}{3})^2}
  \end{pmatrix}
  \tx{,}
\end{align}
and where~$G^*$ is the non-holomorphic Eichler integral of~$G$ normalized as in~\cite{imamoglu-raum-richter-2014,ortiz-raum-richter-2026-preprint}.
Theorem~3.6 of Zwegers's~\cite{zwegers-2001} implies that
\begin{gather}
  \label{eq:definition-of-F}
  F \defeq F^+ + F^-\in\bbM_{\frac{1}{2}}(\rho_3)
  \tx{,}
\end{gather}
where~$\rho_3$ is the type determined by
\begin{gather*}
  \rhoThree(S)
  =
  \zeta_8^{-1}\,
  \begin{psmatrix}
    0 & 1 & 0  \\
    1 & 0 & 0  \\
    0 & 0 & -1
  \end{psmatrix}
  \tx{,}\qquad
  \rhoThree(T)
  =
  \begin{psmatrix}
    \zeta_{24}^{-1} & 0       & 0       \\
    0               & 0       & \zeta_3 \\
    0               & \zeta_3 & 0
  \end{psmatrix}
  \tx{.}
\end{gather*}
Moreover,~\cite{zwegers-2001} provides the transformation laws of~$G$, which implies that~$G\in\rmS_{\frac{3}{2}}(\ov{\rho_3})$.

\subsection{Rankin--Cohen bracket and holomorphic projection}%
\label{ssec:preliminaries:rankin_cohen_projection}

From~\cite{ortiz-raum-richter-2026-preprint} recall the tensor Rankin--Cohen bracket for two vector-valued functions~$F$ and~$G$:
\begin{gather*}
  [F,G]_\nu^\otimes
  \defeq
  \sum_{\mu=0}^\nu
  (-1)^\mu
  \mbinom{k+\nu-1}{\nu-\mu}
  \mbinom{l+\nu-1}{\mu}
  F^{(\mu)}
  \otimes
  G^{(\nu-\mu)}
  \tx{,}
\end{gather*}
where~$F^{(\mu)} \defeq (2\pi i)^{-\mu}\frac{\partial^\mu}{\partial \tau^\mu} F$ is the~$\mu$-th normalized derivative applied elementwise to~$F$.
In particular, if~$F \in \rmM_k(\rho_F)$ and~$G \in \rmM_l(\rho_G)$, then~$[F,G]_\nu^\otimes\in \rmM_{k+l+2\nu}(\rho_F \otimes \rho_G)$.

For a finite-dimensional vector space~$V$ and a continuous function~$F\defcol \HS \ra V$, we write~$\pi_{k}^{\hol}(F)$ for the holomorphic projection of~$F$ in weight~$k$ as defined in~\cite{imamoglu-raum-richter-2014, ortiz-raum-richter-2026-preprint}.
Proposition~4 of~\cite{imamoglu-raum-richter-2014} establishes that~$\pi_{k}^{\hol}(F) = F$ whenever~$F$ is holomorphic.
Furthermore, if~$F$ satisfies appropriate growth conditions and is invariant under the slash operator of weight~$k$ and type~$\rho_F$, then~$\pi_{k}^\hol(F)\in\rmM_k(\rho_F)$.
Consequently, if~$F$ and~$G$ are modular of weights~$k$ and~$l$ and types~$\rho_F$ and~$\rho_G$, respectively, and satisfy analogous growth conditions, it follows that~$\pi_{k+l+2\nu}^{\hol}([F,G]_\nu^\otimes) \in \rmM_{k+l+2\nu}(\rho_F \otimes \rho_G)$.

\section{Vector-Valued Modular Forms}%
\label{sec:vector_valued_modular_forms}
\subsection{Decomposition of a tensor product}%
\label{ssec:vector_valued_modular_forms:tensor_product}

We encounter modular forms of type~$\rho_3\otimes\ov{\rho_3}$, which has the following decomposition.
\begin{lemma}[\cite{imamoglu-raum-richter-2014}, Lemma 1]%
\label{la:vector_valued_modular_forms:rho_three_decomposition}
The representation~$\rhoThree\otimes\ov{\rhoThree}$ is isomorphic to
\begin{gather*}
  \wtd{\sigma}_1 \oplus \wtd{\sigma}_2 \oplus \wtd{\sigma}_6
  \tx{,}
\end{gather*}
where~$\wtd{\sigma}_1$ is the trivial representation and~$\wtd{\sigma}_2$ and~$\wtd{\sigma}_6$ are irreducible representations of dimensions~$2$ and~$6$, respectively.
The summands are realized by the subspaces spanned by the columns of
\begin{gather*}
  \begin{psmatrix}
    1 \\ 0 \\ 0 \\ 0 \\ 1 \\ 0 \\ 0 \\ 0 \\ 1
  \end{psmatrix}
  \tx{,}\quad
  \begin{psmatrix}
    0  & 2  \\
    0  & 0  \\
    0  & 0  \\
    0  & 0  \\
    1  & -1 \\
    0  & 0  \\
    0  & 0  \\
    0  & 0  \\
    -1 & -1
  \end{psmatrix}
  \tx{,}\quad
  \begin{psmatrix}
    0  & 0 & 0 & 0  & 0  & 0 \\
    1  & 1 & 0 & 0  & 0  & 0 \\
    -1 & 1 & 0 & 0  & 0  & 0 \\
    0  & 0 & 1 & 1  & 0  & 0 \\
    0  & 0 & 0 & 0  & 0  & 0 \\
    0  & 0 & 0 & 0  & 1  & 1 \\
    0  & 0 & 1 & -1 & 0  & 0 \\
    0  & 0 & 0 & 0  & -1 & 1 \\
    0  & 0 & 0 & 0  & 0  & 0
  \end{psmatrix}
  \tx{,}
\end{gather*}
respectively.
\end{lemma}

\subsection{Dimensions of spaces of modular forms}%
\label{ssec:vector_valued_modular_forms:dimensions_of_modular_form_spaces}

We are interested in the dimensions of vector-valued cusp forms of weight~$k>2$.
By \zcref{la:vector_valued_modular_forms:rho_three_decomposition}, we have the decomposition
\begin{gather*}
  \rmM_{2+2\nu}\bigl( \rhoThree\otimes\ov{\rhoThree} \bigr)
  \cong
  \rmM_{2+2\nu}(\wtd\sigma_1)
  \oplus
  \rmM_{2+2\nu}(\wtd\sigma_2)
  \oplus
  \rmM_{2+2\nu}(\wtd\sigma_6)
  \tx{.}
\end{gather*}
The dimensions of the summands are calculated using the dimension formula on p.~228 of~\cite{borcherds-1999} and the Eisenstein subspace characterization from~\cite{westerholt-raum-2017}.
For brevity, the details are omitted, and we record only the following proposition.
\begin{proposition}%
\label{prop:vector_valued_modular_forms:target_space_dimensions}
Let~$k > 2$ be even.
For~$\wtd\sigma_1$, we have
\begin{gather*}
  \dim \rmM_k(\wtd\sigma_1)
  =
  \begin{cases}
    \bigl\lfloor \frac{k}{12} \bigr\rfloor\tx{,}
    \\
    \bigl\lfloor \frac{k}{12} \bigr\rfloor + 1\tx{,}
  \end{cases}
  \quad\tx{and}\quad
  \dim \rmS_k(\wtd\sigma_1)
  =
  \begin{cases}
    \bigl\lfloor \frac{k}{12} \bigr\rfloor - 1\tx{,} &
    \qquad
    \tx{if\/ } k \equiv 2 \quad\pmod{12}\tx{;}
    \\
    \bigl\lfloor \frac{k}{12} \bigr\rfloor\tx{,} &
    \qquad
    \tx{if\/ } k \not\equiv 2 \quad\pmod{12}\tx{.}
  \end{cases}
\end{gather*}
For~$\wtd\sigma_2$, we have
\begin{gather*}
  \dim \rmM_k(\wtd\sigma_2)
  =
  \bigl\lfloor
  \mfrac{k+4}{6}
  \bigr\rfloor
  \quad\tx{and}\quad
  \dim \rmS_k(\wtd\sigma_2)
  =
  \bigl\lfloor
  \mfrac{k+4}{6}
  \bigr\rfloor
  -
  1
  \tx{.}
\end{gather*}
For~$\wtd\sigma_6$, we have
\begin{gather*}
  \dim \rmM_k(\wtd\sigma_6)
  =
  \mfrac{k}{2}
  \quad\tx{and}\quad
  \dim \rmS_k(\wtd\sigma_6)
  =
  \mfrac{k}{2}
  -
  1
  \tx{.}
\end{gather*}
In particular,
\begin{gather*}
  \rmS_4(\wtd\sigma_1)
  =
  \rmS_6(\wtd\sigma_1)
  =
  \rmS_4(\wtd\sigma_2)
  =
  \rmS_6(\wtd\sigma_2)
  =
  0
  \tx{.}
\end{gather*}
\end{proposition}

\section{Fourier Coefficients of Holomorphic Projections}%
\label{sec:holomorphic_projection_fourier_coefficients}

Recalling the functions $F$ and $G$ from \zcref{ssec:preliminaries:mock_theta_vector}, we apply the vector-valued holomorphic projection operator to obtain
\begin{gather*}
  \pi_{2+2\nu}^{\hol}
  \bigl(
  [F,G]_\nu^\otimes
  \bigr)
  \in
  \rmM_{2+2\nu}\bigl( \rhoThree\otimes\ov{\rhoThree} \bigr)
  \tx{.}
\end{gather*}
The holomorphic part yields a direct contribution via $[F^+,G]_\nu^\otimes$, whereas the nonholomorphic part is determined by applying the projection formulas from~\cite{ortiz-raum-richter-2026-preprint}.

Throughout, let~$\frake_1,\frake_2,\frake_3$ be the standard basis of the representation space of~$\rhoThree$.

\subsection{Coefficients of holomorphic terms}%
\label{ssec:holomorphic_projection_fourier_coefficients:holomorphic}

We first determine the Fourier coefficients of~$[F^+,G]_\nu^\otimes$ in terms of the coefficients of the mock theta functions.
Define the coefficients~$c_{r,s,\nu}^+(N)$ by
\begin{gather*}
  [F^+,G]_\nu^\otimes
  =
  \sum_{r,s=1}^3
  \sum_N
  c_{r,s,\nu}^+(N)\,
  e(N\tau)\,
  \frake_r\otimes\frake_s
  \tx{.}
\end{gather*}

Write
\begin{gather*}
  F^+(\tau)
  =
  \sum_{r=1}^3
  \sum_{\td{m}}
  c^+\bigl( F^{[r]};\, \td{m} \bigr)\,
  e(\td{m}\tau)\,
  \frake_r
  \quad\tx{and}\quad
  G(\tau)
  =
  \sum_{s=1}^3
  \sum_{\td{n}>0}
  c\bigl( G^{[s]};\, \td{n} \bigr)\,
  e(\td{n}\tau)\,
  \frake_s
  \tx{.}
\end{gather*}
The nonzero coefficients of~$F^+$ are related to~$f$ and~$\omega$ by
\begin{gather*}
  c^+\bigl(F^{[1]};n-\tfrac{1}{24}\bigr)
  =
  c_f(n)
  \tx{,}\qquad
  c^+\bigl(F^{[2]};\tfrac{n}{2}+\tfrac{1}{3}\bigr)
  =
  2c_\omega(n)
  \tx{,}\qquad
  c^+\bigl(F^{[3]};\tfrac{n}{2}+\tfrac{1}{3}\bigr)
  =
  2(-1)^n\,
  c_\omega(n)
  \tx{.}
\end{gather*}
Setting~$h_1=1$ and~$h_2=h_3=2$, we define for~$m\in 6\ZZ+h_r$,
\begin{gather*}
  \epsilon_1(m)
  \defeq
  -1
  \tx{,}\qquad
  \epsilon_2(m)
  \defeq
  (-1)^{\frac{m-2}{6}}
  \tx{,}\qquad
  \epsilon_3(m)
  \defeq
  -1
  \tx{.}
\end{gather*}
The coefficients of the theta series in~\eqref{eq:definition-of-G} are then given by
\begin{gather}
  \label{eq:holomorphic_projection:theta_coefficients}
  c\bigl( G^{[r]};\, \td{m} \bigr)
  =
  \begin{cases}
    \epsilon_r(m)\,\frac{m}{6}
    \tx{,}
     &
    \tx{if }\td{m}=\frac{m^2}{24}\tx{ with }m\in 6\ZZ+h_r
    \tx{;}
    \\
    0
    \tx{,}
     &
    \tx{otherwise}
    \tx{.}
  \end{cases}
\end{gather}

In connection with the Rankin--Cohen bracket, we consider the polynomial
\begin{gather*}
  P_\nu(X,Y)
  \defeq
  \sum_{\mu=0}^\nu
  (-1)^\mu\,
  \mbinom{\nu-\frac{1}{2}}{\nu-\mu}
  \mbinom{\nu+\frac{1}{2}}{\mu}\,
  X^\mu
  Y^{\nu-\mu}
  \tx{.}
\end{gather*}
The definition of the Rankin--Cohen bracket implies that
\begin{gather}
  \label{eq:holomorphic_projection:holomorphic_component_coefficients}
  c_{r,s,\nu}^+(N)
  =
  \sum_{\substack{\td{m},\;\td{n}>0\\ \td{m}+\td{n}=N}}
  \mspace{-8mu}
  c^+\bigl( F^{[r]};\, \td{m} \bigr)\,
  c\bigl( G^{[s]};\, \td{n} \bigr)\,
  P_\nu(\td{m},\td{n})
  \tx{.}
\end{gather}
Substituting \eqref{eq:holomorphic_projection:theta_coefficients} into \eqref{eq:holomorphic_projection:holomorphic_component_coefficients}, we obtain
\begin{gather}
  \label{eq:holomorphic_projection:holomorphic_component_coefficients_theta_inserted}
  c_{r,s,\nu}^+(N)
  =
  \mfrac{1}{6}
  \sum_{m\in 6\ZZ+h_s}
  \epsilon_s(m)\,
  m\,
  c^+\bigl(F^{[r]};\, N-\mfrac{m^2}{24}\bigr)\,
  P_\nu\bigl(N-\mfrac{m^2}{24},\mfrac{m^2}{24}\bigr)
  \tx{.}
\end{gather}
Observe the sum is finite since~$c_f(n)=c_\omega(n)=0$ for~$n\not\in\ZZ_{\ge0}$.

Furthermore, we evaluate the constant coefficient of~$[F^+,G]_\nu^\otimes$ explicitly in order to identify the vector-valued Eisenstein series in \zcref{sec:proof_main_theorem}.
For~$N=0$, the only nonzero term in \eqref{eq:holomorphic_projection:holomorphic_component_coefficients_theta_inserted} arises when $r=s=1$ and $m=1$.
Since~$c^+(F^{[1]};-\frac{1}{24})=c_f(0)=1$ and~$\epsilon_1(1)=-1$, applying Vandermonde's identity to $P_\nu\bigl(-\frac{1}{24},\frac{1}{24}\bigr)$ yields
\begin{gather}
  \label{eq:holomorphic_projection:holomorphic_constant_term}
  c_{r,s,\nu}^+(0)
  =
  \begin{cases}
    -\mfrac{1}{6\cdot 24^\nu}
    \mbinom{2\nu}{\nu}
    \tx{,}
     &
    \tx{if } (r,s)=(1,1)
    \tx{;}
    \\
    0
    \tx{,}
     &
    \tx{if } (r,s)\ne(1,1)
    \tx{.}
  \end{cases}
\end{gather}

\subsection{Coefficients of nonholomorphic terms}%
\label{ssec:holomorphic_projection_fourier_coefficients:nonholomorphic}

We next determine the Fourier coefficients of the holomorphic projection of~$[F^-,G]_\nu^\otimes$.
We define the coefficients~$c_{r,s,\nu}^-(N)$ by
\begin{gather*}
  \pi_{2+2\nu}^{\hol}
  \bigl(
  [F^-,G]_\nu^\otimes
  \bigr)
  =
  \sum_{r,s=1}^3
  \sum_{N>0}
  c_{r,s,\nu}^-(N)\,
  e(N\tau)\,
  \frake_r\otimes\frake_s
  \tx{.}
\end{gather*}

\begin{remark}%
\label{rm:holomorphic_projection:nonholomorphic_constant_term}
In contrast to~\cite{ortiz-raum-richter-2026-preprint}, where a pure power contribution of~$y^{-\kappa}$ arises in the nonholomorphic part due to~$1-k<0$, no such term appears in our setting since~$1-k>0$.
\end{remark}

We apply Proposition 3.7 of~\cite{ortiz-raum-richter-2026-preprint} with~$k=\frac{1}{2}$ and~$l=\frac{3}{2}$, along with the binomial identity
\begin{gather*}
  \mbinom{\nu-\frac{1}{2}}{\nu-t}
  \mbinom{\nu+\frac{1}{2}}{t}
  =
  \mbinom{\nu-\frac{1}{2}}{\nu}
  \mbinom{2\nu+1}{2t}
\end{gather*}
to obtain
\begin{gather}
  \label{eq:holomorphic_projection:nonholomorphic_specialized}
  \pi_{2+2\nu}^{\hol}
  \bigl(
  \bigl[
    \Ga\bigl(\tfrac{1}{2},4\pi |n| y \bigr) e(n\tau),\,
    e(\td{n}\tau)
  \bigr]_\nu
  \bigr)
  =
  \sqrt{\pi}
  \mbinom{\nu-\frac{1}{2}}{\nu}\,
  \mfrac{
    \bigl(
    \sqrt{\td{n}}
    -
    \sqrt{|n|}
    \bigr)^{2\nu+1}
  }{
    \sqrt{\td{n}}
  }\,
  e\bigl((n+\td{n})\tau\bigr)
  \tx{,}
\end{gather}
where~$\Ga$ denotes the upper incomplete Gamma function.

Recalling that~$F^-=-2\sqrt{6} G^*(\tau)$, the definition of the Eichler integral implies that
\begin{gather*}
  F^-(\tau)
  =
  (4\pi)^{-\frac{1}{2}}\,
  2\sqrt{6}\,
  \sum_{r=1}^3
  \sum_{n<0}
  c\bigl( G^{[r]};\, |n| \bigr)\,
  |n|^{-\frac{1}{2}}\,
  \Ga\bigl( \tfrac{1}{2},4\pi |n|y \bigr)
  e(n\tau)\,
  \frake_r
  \tx{.}
\end{gather*}
For~$N>0$, equation~\eqref{eq:holomorphic_projection:nonholomorphic_specialized} gives
\begin{gather}
  \label{eq:holomorphic_projection:nonholomorphic_component_coefficients}
  c_{r,s,\nu}^-(N)
  =
  \sqrt{6}
  \mbinom{\nu-\frac{1}{2}}{\nu}
  \sum_{\substack{n<0,\;\td{n}>0\\ n+\td{n}=N}}
  c\bigl( G^{[r]};\, |n| \bigr)\,
  c\bigl( G^{[s]};\, \td{n} \bigr)\,
  \mfrac{
    \bigl(
    \sqrt{\td{n}}
    -
    \sqrt{|n|}
    \bigr)^{2\nu+1}
  }{
    \sqrt{|n|\td{n}}
  }
  \tx{,}
\end{gather}
and substituting~\eqref{eq:holomorphic_projection:theta_coefficients} into~\eqref{eq:holomorphic_projection:nonholomorphic_component_coefficients}, setting~$|n|=\frac{a^2}{24}$ and~$\td{n}=\frac{b^2}{24}$, and simplifying yields
\begin{gather}
  \label{eq:holomorphic_projection:nonholomorphic_component_coefficients_theta_inserted}
  \begin{aligned}
  c_{r,s,\nu}^-(N)
  ={} &
  \mfrac{1}{3\cdot 24^\nu}
  \mbinom{\nu-\frac{1}{2}}{\nu}
  \sum_{\substack{
      a\in 6\ZZ+h_r,\; b\in 6\ZZ+h_s\\
      b^2-a^2=24N
    }}
  \mspace{-18mu}
  \epsilon_r(a)
  \epsilon_s(b)\,
  \sgn(ab)\,
  \bigl(
  |b|-|a|
  \bigr)^{2\nu+1}
  \tx{.}
  \end{aligned}
\end{gather}
To write \eqref{eq:holomorphic_projection:nonholomorphic_component_coefficients_theta_inserted} as a divisor sum, we change variables by setting~$d=b-a$ and~$e=b+a$, so that
\begin{gather*}
  de=24N
  \quad\tx{and}\quad
  a=\mfrac{e-d}{2}
  \quad\tx{and}\quad
  b=\mfrac{e+d}{2}
  \tx{.}
\end{gather*}
Since~$N>0$, the two divisors~$d$ and~$e$ share the same sign, and hence
\begin{gather*}
  |b|-|a|
  =
  \min\bigl(|d|,|e|\bigr)
  \tx{.}
\end{gather*}
In addition,
\begin{gather*}
  \sgn(ab)
  =
  \sgn(e^2-d^2)
  \tx{.}
\end{gather*}
Altogether, \eqref{eq:holomorphic_projection:nonholomorphic_component_coefficients_theta_inserted} becomes
\begin{gather}
  \label{eq:holomorphic_projection:nonholomorphic_component_coefficients_divisor_sum}
  \begin{aligned}
  c_{r,s,\nu}^-(N)
  ={} &
  \mfrac{1}{3\cdot 24^\nu}
  \mbinom{\nu-\frac{1}{2}}{\nu}
  \sum_{\substack{
      d,e\in\ZZ\\
      de=24N\\
      \frac{e-d}{2}\in 6\ZZ+h_r\\
      \frac{e+d}{2}\in 6\ZZ+h_s
    }}
  \epsilon_r\bigl(\mfrac{e-d}{2}\bigr)
  \epsilon_s\bigl(\mfrac{e+d}{2}\bigr)\;
  \sgn(e^2-d^2)\,
  \min\bigl(|d|,|e|\bigr)^{2\nu+1}
  \tx{,}
  \end{aligned}
\end{gather}
where the sum is empty (and the coefficient vanishes) if~$24N$ is not an integer.

\subsection{Coefficients of the diagonal components}%
\label{ssec:holomorphic_projection_fourier_coefficients:diagonal_coefficients}

We now evaluate the coefficients of the diagonal components.
By \zcref{la:vector_valued_modular_forms:rho_three_decomposition}, the coefficients of the mixed components~$\frake_i\otimes\frake_j$ with~$i\ne j$ receive contributions from cusp forms when~$k>2$; consequently, we exclude them from our analysis.
The components~$\frake_1\otimes\frake_1$ and~$\frake_2\otimes\frake_2$, give rise to the recursions in \zcref{mainthm:mock_theta_recursions} for the coefficients of~$f$ and~$\omega$, respectively.
Since~$T$ exchanges the second and third components, the component~$\frake_3\otimes\frake_3$ yields results identical to those of~$\frake_2\otimes\frake_2$, and hence it is also omitted.

Explicitly, from~\eqref{eq:holomorphic_projection:holomorphic_component_coefficients_theta_inserted}, the contribution of the holomorphic part equals
\begin{gather}
  \label{eq:holomorphic_projection:diagonal_holomorphic_components}
  \begin{alignedat}{2}
  c_{1,1,\nu}^+(N)
   & =
  -\mfrac{1}{6}
  \sum_{m\in6\ZZ+1}
  m\,
  c_f\bigl(N+\mfrac{1-m^2}{24}\bigr)\,
  P_\nu\bigl(N-\mfrac{m^2}{24},\mfrac{m^2}{24}\bigr)
  \quad
  &&
  \tx{for } N\in \ZZ_{>0}
  \tx{;}
  \\
  c_{2,2,\nu}^+(N)
   & =
  \mfrac{1}{3}
  \sum_{m\in6\ZZ+2}
  (-1)^{\frac{m-2}{6}}\,
  m\,
  c_\omega\bigl(2N-\mfrac{m^2}{12}-\mfrac{2}{3}\bigr)\,
  P_\nu\bigl(N-\mfrac{m^2}{24},\mfrac{m^2}{24}\bigr)
  \quad
  &&
  \tx{for } N\in\tfrac{1}{2}\ZZ_{>0}
  \tx{.}
  \end{alignedat}
\end{gather}
From~\eqref{eq:holomorphic_projection:nonholomorphic_component_coefficients_divisor_sum}, the nonholomorphic part evaluates to
\begin{gather}
  \label{eq:holomorphic_projection:diagonal_nonholomorphic_components}
  c_{1,1,\nu}^-(N)
  =
  \mfrac{1}{3\cdot24^\nu}
  \mbinom{\nu-\frac{1}{2}}{\nu}\,
  \lambda^f_{2\nu+1}(N)
  \quad\tx{and}\quad
  c_{2,2,\nu}^-(N)
  =
  \mfrac{1}{3\cdot24^\nu}
  \mbinom{\nu-\frac{1}{2}}{\nu}\,
  \lambda^\omega_{2\nu+1}(N)
  \tx{,}
\end{gather}
where $\lambda^f_\nu(N)$ and $\lambda^\omega_\nu(N)$ are as in \zcref{mainthm:mock_theta_recursions}.

\section{Proof of \zcref{mainthm:mock_theta_recursions}}%
\label{sec:proof_main_theorem}

Applying holomorphic projection to the de\-com\-po\-si\-tion~$[F,G]_\nu^\otimes = [F^+,G]_\nu^\otimes + [F^-,G]_\nu^\otimes$ for~$F$ and~$G$ as in \zcref{ssec:preliminaries:mock_theta_vector} gives
\begin{gather}
  \label{eq:proof_main_theorem:projected_bracket_identity}
  [F^+,G]_\nu^\otimes
  =
  -
  \pi_{2+2\nu}^{\hol} \bigl(
  [F^-,G]_\nu^\otimes
  \bigr)
  +
  \pi_{2+2\nu}^{\hol} \bigl(
  [F,G]_\nu^\otimes
  \bigr)
  \tx{.}
\end{gather}
The left hand side contains the desired weighted sums of mock theta coefficients evaluated in \zcref{ssec:holomorphic_projection_fourier_coefficients:holomorphic}, while the first term on the right hand side was evaluated in \zcref{ssec:holomorphic_projection_fourier_coefficients:nonholomorphic}.
Therefore, it remains only to identify the modular form
\begin{gather*}
  \pi_{2+2\nu}^{\hol}\bigl(
  [F,G]_\nu^\otimes
  \bigr)
  \in
  \rmM_{2+2\nu}\bigl( \rhoThree\otimes\ov{\rhoThree} \bigr)
  \tx{.}
\end{gather*}
The relations in \zcref{mainthm:mock_theta_recursions} arise from the diagonal components of~$\rhoThree\otimes\ov{\rhoThree}$, which span the subrepresentation~$\wtd{\sigma}\defeq\wtd{\sigma}_1\oplus\wtd{\sigma}_2$.
Using the basis~$\frake_i$ of~$\rho_3$ from \zcref{sec:holomorphic_projection_fourier_coefficients}, we define the basis elements
\begin{gather*}
  \frakd_1 \defeq \frake_1\otimes\frake_1
  \tx{,}\quad
  \frakd_2 \defeq \frake_2\otimes\frake_2
  \tx{,}\quad
  \frakd_3 \defeq \frake_3\otimes\frake_3
  \tx{.}
\end{gather*}
Since \zcref{prop:vector_valued_modular_forms:target_space_dimensions} asserts that~$\rmS_4(\wtd{\sigma}) = \rmS_6(\wtd{\sigma}) = 0$, it suffices to match the constant coefficients of the diagonal components with vector-valued Eisenstein series that span~$\rmM_{2 + 2\nu}(\wtd\sigma)$.

For even~$k>2$, let~$E_k$ denote the normalized scalar Eisenstein series on~$\SL{2}(\ZZ)$:
\begin{gather*}
  E_4(\tau)
  =
  1
  +
  240 \sum_{n=1}^\infty \sigma_3(n)\, e(n \tau)
  \tx{,}\qquad
  E_6(\tau)
  =
  1
  -
  504 \sum_{n=0}^\infty \sigma_5(n)\, e(n \tau)
  \tx{.}
\end{gather*}
With respect to the basis~$\frakd_1,\frakd_2,\frakd_3$, the representation~$\wtd{\sigma}$ acts as a permutation, given explicitly by
\begin{gather*}
  \wtd{\sigma}(S)
  =
  \begin{psmatrix}
    0 & 1 & 0 \\
    1 & 0 & 0 \\
    0 & 0 & 1
  \end{psmatrix}
  \tx{,}\quad
  \wtd{\sigma}(T)
  =
  \begin{psmatrix}
    1 & 0 & 0 \\
    0 & 0 & 1 \\
    0 & 1 & 0
  \end{psmatrix}
  \tx{.}
\end{gather*}
Note that the first two components are exchanged by~$S$, the last component is fixed by~$S$, and~$T$ fixes the first component while exchanging the remaining two.
We conclude that~$\wtd{\sigma}$ is the induction of the trivial character from the congruence subgroup of level~$2$ to~$\SL{2}(\ZZ)$ (see also Section 2.3 of~\cite{ortiz-raum-richter-2026-preprint}).
Specifically, the first components of modular forms of type~$\wtd\sigma$ are in one-to-one correspondence with scalar-valued modular forms for~$\Ga_0(2)$ of the same weight.
The second and third components of such a modular form arise from the first one by letting~$S$ and~$S T$ act.
Thus we have an Eisenstein series~$E_{k,\frakd_1}$ of type~$\wtd\sigma$ with components
\begin{gather}
  \label{eq:proof_main_theorem:diagonal_eisenstein_components}
  E_{k,\frakd_1}^{[1]}\!(\tau)
  \defeq
  \mfrac{2^k E_k(2\tau)-E_k(\tau)}{2^k-1}
  \tx{,}\quad
  E_{k,\frakd_1}^{[2]}\!(\tau)
  \defeq
  \mfrac{E_k(\tau/2)-E_k(\tau)}{2^k-1}
  \tx{,}\quad
  E_{k,\frakd_1}^{[3]}\!(\tau)
  \defeq
  \mfrac{E_k((\tau+1)/2)-E_k(\tau)}{2^k-1}
  \tx{,}
\end{gather}
whose constant coefficient equals~$\frakd_1$.
We let~$c(E_{k,\frakd_1}^{[i]};N)$ denote the~$N$-th Fourier coefficient of the~$i$\thdash{} component for~$i \in \{1,2,3\}$.
Inserting the classical weight~$4$ and~$6$ Eisenstein coefficients into~\eqref{eq:proof_main_theorem:diagonal_eisenstein_components} yields
\begin{gather}
\label{eq:proof_main_theorem:diagonal_eisenstein_fourier_coefficients}
\begin{alignedat}{2}
  c\bigl(E_{4,\frakd_1}^{[1]};\, N\bigr)
  & =
  -
  16
  \bigl(
  \sigma_3(N)
  -
  16\,\sigma_3\bigl( \tfrac{N}{2} \bigr)
  \bigr)
  \tx{,}
  \qquad
  &
  c\bigl(E_{4,\frakd_1}^{[2]};N\bigr)
   & =
  16
  \bigl(
  \sigma_3(2N)
  -
  \sigma_3(N)
  \bigr)
  \tx{,}
  \\
  c\bigl(E_{6,\frakd_1}^{[1]};N\bigr)
   & =
  8
  \bigl(
  \sigma_5(N)
  -
  64\, \sigma_5\bigl( \tfrac{N}{2} \bigr)
  \bigr)
  \tx{,}
  \qquad
  &
  c\bigl(E_{6,\frakd_1}^{[2]};N\bigr)
   & =
  -8
  \bigl(
  \sigma_5(2N)
  -
  \sigma_5(N)
  \bigr)
  \tx{.}
\end{alignedat}
\end{gather}

Recall from~\eqref{eq:holomorphic_projection:holomorphic_constant_term} (see also \zcref{rm:holomorphic_projection:nonholomorphic_constant_term}) that the constant Fourier coefficient of~$\pi_{2+2\nu}^{\hol}([F,G]_\nu^\otimes)$ is given by~$\frac{-1}{6\cdot 24^\nu} \binom{2\nu}{\nu}\, \frakd_1$.
Thus, the diagonal components of~$\pi_{2+2\nu}^{\hol}([F,G]_\nu^\otimes)$ are
\begin{gather}
  -\mfrac{1}{6\cdot 24^\nu}
  \mbinom{2\nu}{\nu}\,
  \bigl(
  E_{2+2\nu,\frakd_1}^{[1]}\!(\tau)\,\frakd_1
  +
  E_{2+2\nu,\frakd_1}^{[2]}\!(\tau)\,\frakd_2
  +
  E_{2+2\nu,\frakd_1}^{[3]}\!(\tau)\,\frakd_3
  \bigr)
  \tx{.}
\end{gather}
This allows us to compare the~$\frake_1\otimes\frake_1$- and~$\frake_2\otimes\frake_2$-components in~\eqref{eq:proof_main_theorem:projected_bracket_identity}.
Recall also the divisor sums~$\lambda_{\nu}^f$ and~$\lambda_{\nu}^\omega$ from \zcref{mainthm:mock_theta_recursions}.
Using \eqref{eq:holomorphic_projection:diagonal_holomorphic_components} and~\eqref{eq:holomorphic_projection:diagonal_nonholomorphic_components}, we find
\begin{gather}
  \label{eq:proof_main_theorem:general_recursions}
  \begin{alignedat}{2}
  c_{1,1,\nu}^+(N)
  & =
  -
  \mfrac{1}{3\cdot24^\nu}
  \mbinom{\nu-\frac{1}{2}}{\nu}\,
  \lambda^f_{2\nu+1}(N)
  -
  \mfrac{1}{6\cdot24^\nu}
  \mbinom{2\nu}{\nu}\,
  c\bigl( E_{2+2\nu,\frakd_1}^{[1]};\, N \bigr)
  \tx{,}
  \quad
  &&
  \tx{for } N\in\ZZ_{>0}
  \tx{;}
  \\
  c_{2,2,\nu}^+(N)
  & =
  -
  \mfrac{1}{3\cdot24^\nu}
  \mbinom{\nu-\frac{1}{2}}{\nu}\,
  \lambda^\omega_{2\nu+1}(N)
  -
  \mfrac{1}{6\cdot24^\nu}
  \mbinom{2\nu}{\nu}\,
  c\bigl( E_{2+2\nu,\frakd_1}^{[2]};\, N \bigr)
  \tx{,}
  \quad
  &&
  \tx{for } N\in\tfrac{1}{2}\ZZ_{>0}
  \tx{.}
  \end{alignedat}
\end{gather}
After simplification, Equations~\eqref{eq:holomorphic_projection:diagonal_nonholomorphic_components}, \eqref{eq:proof_main_theorem:diagonal_eisenstein_fourier_coefficients}, and~\eqref{eq:proof_main_theorem:general_recursions} yield the four recursions stated in \zcref{mainthm:mock_theta_recursions}.

\begin{appendix}

\newenvironment{prompt}
  {\ttfamily\small
   \spaceskip=\fontdimen2\font plus .6\fontdimen2\font minus .3\fontdimen2\font
   \emergencystretch=1.5em}
  {}

\section{AI-Assisted Research}%
\label{sec:ai-contribution}

This appendix records how large language models (AI) contributed to the present work.
We provide enough detail to make the workflow transparent and, most of all, reproducible.

Our work sets itself apart from two lines of contemporary research.
The first couples language models with proof assistants such as \texttt{Lean} or \texttt{Isabelle}, so that every inference is machine-checked~\cite{ono-et-al-2026a-preprint,ono-et-al-2026b-preprint,google-deepmind-2026-preprint}.
This guarantees formal correctness, at the price of confining the argument to a formalizable fragment and of a considerable formalization effort.
We make no use of formal verification: our proofs are written in ordinary mathematical prose and were verified by the human authors.
The second line trains or deploys specialized, frequently proprietary systems of substantial cost dedicated to mathematical problem solving~\cite{patel-2026-preprint,alon-bloom-gowers-litt-sawin-shankar-tsimerman-wang-wood-2026-preprint,openai-2026}.
We instead rely on a general-purpose, off-the-shelf model accessed through an ordinary consumer subscription; tests with open-weight models reproduced the complete argument as well, albeit in a more technical formulation.
This unlocks a workflow in which affordable and widely available tools contribute to research mathematics under human supervision.
The mathematical craft lies within a handful of inspectable instructions that the reader can examine, adapt, and iterate on.

\paragraph{Chat interfaces and agent harnesses}

Most readers will have met a language model through a web chat interface: one types a request in natural language and receives a textual reply.
In these cases the model sees nothing beyond the text of the conversation, with the exception of the by-now common \enquote{memory} features.
A \emph{harness} is a program that wraps the model and equips it with tools to read and write files in a working directory, to run code, or to query the web.
It feeds the result of each tool back to the model, so that the model can act on more than the information it was trained on.
The model thereby becomes an \emph{agent} acting within the project.
For example, one of the prompts in \zcref{ssec:ai-contribution:prompts} instructs the model to amend~\texttt{AGENTS.md}, a plain-text description of the manuscript that the harness presents to the agent at the start of every session.
A chat interface would answer with advice; the agent instead edits the file.
An agent operates in a loop: it is given a task, chooses an action, the harness executes it and returns the outcome, and the cycle repeats until the task is deemed complete.
A single session for the present work typically ran through many such steps over a long interaction.

\paragraph{Agent skills}

A \emph{skill} is a self-contained set of instructions that the agent loads on demand when a task calls for it.
Skills let us fix, once and explicitly, how each recurring mathematical task is to be carried out, rather than entrusting it to the model's defaults.
The three skills written for this work encode three competences of a working mathematician: reading and comparing the references that enter the argument, testing formulas and simplifications by computation, and typesetting clean mathematics whose structure is laid out clearly in the source.
\zcref{ssec:ai-contribution:skills} describes each.
This agent-and-skill methodology is by now standard in AI-assisted software development, from which we have borrowed it; its transfer to mathematical research is what the present work explores.

\paragraph{Setup}

We used the model GPT-5.5 at \enquote{xhigh} reasoning effort, accessed through the Codex harness via a standard subscription; other state-of-the-art models and open harnesses are equally realistic choices.
As external material we supplied the two source papers~\cite{imamoglu-raum-richter-2014,ortiz-raum-richter-2026-preprint} together with the references of Borcherds~\cite{borcherds-1999} and of Raum~\cite{westerholt-raum-2017} that furnish the relevant dimension formulas.
The mathematics was developed in~12 prompts, reproduced in full in \zcref{ssec:ai-contribution:prompts}.

\subsection{Supplementary Material}%
\label{ssec:ai-contribution:supplementary-material}

The supplementary material consists of
\begin{enumerateroman*}
\item the three agent skills of \zcref{ssec:ai-contribution:skills} in full,
\item the initial project files, that is, the customary LaTeX templates, including the shortcut definitions, routinely used by one of the authors,
\item the diffs that track the evolution of the accompanying AI notes, and
\item the session logs recording the interaction with the agent harness.
\end{enumerateroman*}
It is provided at the second author's homepage, at~\url{martin.raum-brothers.eu/software}.

\subsection{Prompts}%
\label{ssec:ai-contribution:prompts}

The development proceeded through the prompts reproduced below, in the order in which they were issued.
In the third and fourth interaction the agent asked a question or requested confirmation before writing the planned changes to file, which accounts for the extra prompt in each.

\paragraph{Combining two papers}

\begin{prompt}
We are planing a mathematical manuscript that combines ideas from the two papers that you find in external\_material. In short, we take the vector-valued Rankin-Cohen bracket (second paper) and apply it to the vector-valued completion of the third order mock theta series (first paper) and its image under the xi operator, which is a theta series. The general structure of this work should follow the second paper. In the preliminaries you can collect the information on the basic constructions. In a section on vector-valued modular forms you can include the representation theoretic computation from the first paper. This will require a section on induced representations similar to the one in the second paper, but with different representations. The calculations of the holomorphic projection, for now, can be cited from the second paper. They go to a separate section. Then finally, we need to apply the Rankin-Cohen brackets, evaluate the holomorphic projection, and simplify. We have to later decide which of the target spaces of modular forms vanish to derive the desired weighted recursions; reserve a placehold for this in both the section on vector-valued modular forms and in the final section, where we assemble the proof.
\end{prompt}

\paragraph{Outline paper in Agent File}

\begin{prompt}
Amend the AGENTS.md with a concise description of the manuscript's goal per its abstract and introduction sketch; Add also a very brief summary of each sections (intended) content.
\end{prompt}

\paragraph{Ingesting Borcherds' dimension formula}

\begin{prompt}
We have changed the idea about section 2 slightly: For the sake of brevity, we skip the part on the idenitfication of irreducible components as genuine inductions. Instead, we want to use a dimension formula by Borcherds to show that specific spaces of modular forms and of cusp forms vanish. Plan to remove or update the affected sections as appropriate.

Go ahead using the citation by Borcherds; Use the citation summary skill to extract the relevant dimension formula from external\_material.
\end{prompt}

\paragraph{Ingesting cusp form dimension formula}

\begin{prompt}
The dimension of spaces of cusp forms lacks justification. To fix this we now ingest the relevant the characterization of Eisenstein series by Raum. The paper is available in external\_material. Use the citation summary skill to plan extracting the characterization of Eisenstein series. Note that the condition \$k \textbackslash{}ge 2\$ is a typo and should read \$k > 2\$. The case of \$k = 2\$ is stated separately.

Implement.
\end{prompt}

\paragraph{Test dimension formula}

\begin{prompt}
Add the dimension formula for sigma\_1 to Proposition 2.2 and use the verification skill to check the whole proposition numerically.
\end{prompt}

\paragraph{Determine contribution from non-holomorphic terms}

\begin{prompt}
We work on Section 3 and need to insert the Fourier coefficient formulas for the contribution of the nonholomorphic part. They are determined in the manuscript on Hurwitz class numbers, available in external\_material. Use the citation summary skill to extract the relevant formula for general l and k, insert it as a citation and then specialize to the values of k and l used in the manuscript. Verify the resulting simplied formula numerically using the verification skill. Note (and remark in the manuscript) that as opposed to the case of Hurwitz class numbers the constant term does not receive any contribution.
\end{prompt}

\paragraph{Ingest modular completion of mock theta series}

\begin{prompt}
We have to expand on Section 1.2 and determine the exact modular completion of the third order mock theta series. Use the citation summary skill to extract it from the third order mock theta paper in external\_material. Then use the verification skill to check numerically modular invariance of the completion by inserting points close to i on the upper half plane into the power series expansion and its image under the slash action of S. Formulate the modular completion in terms of the non-holomorphic Eichler integral of the theta series, which you find in both the preliminaries of the mock theta and the Hurwitz class numbers paper in external\_material (use again the citation summary skill and the verification skill).
\end{prompt}

\paragraph{Simplify contribution from non-holomorphic terms}

\begin{prompt}
Insert the values for the coefficients of the theta series into the contribution (3.3) of the non-holomorphic terms and then simplify the sum. By writing N as a product of a difference and a sum, we effectively expect a divisor sum, similar to the ones in the mock theta and the hurwitz papers, which are available in external\_material.
\end{prompt}

\paragraph{Determine contribution from holomorphic terms}

\begin{prompt}
We further resolve Section 3.2: Add a section before the one on contributions from non-holomorphic terms that summarizes and simplifies the contributions from the holomorphic terms.
\end{prompt}

\paragraph{Determine contribution from Eisenstein series}

\begin{prompt}
From the holomorphic contributions in Section 3.2 determine the contribution to the constant term. Then project to the span of sigma1 and sigma2 and match the corresponding Eisenstein series, which are similar to the ones in the previous paper on mock theta functions in external\_material. Note however that we are in weight 2 + 2nu > 2, so that no quasi modular Eisenstein series are needed.
\end{prompt}

\paragraph{Assemble recursions}

\begin{prompt}
We now finish the paper: Discard Section 3.4 and implement it. In Section 4 we combine the contributions determined in Section 3 for nu = 1 and nu = 2. In both cases the modular forms contributions from sigma1 and sigma2 are Eisenstein series, since the corresponding space of cusp forms vanish. Deduce the resulting recursion relations for the coefficients of the mock theta functions f and omega. When done, write a main theorem in the introduction that clearly states this recursion without forward references to the body of the paper. Follow the style of the main theorem of the paper on Hurwitz class numbers, which you find in external\_material.
\end{prompt}

\paragraph{Test recursions}

\begin{prompt}
Use the verification skill to implement the recursions of the main theorem in Julia and test them for many values.
\end{prompt}

\subsection{Agent Skills}%
\label{ssec:ai-contribution:skills}

Each skill is a substantial document of potentially hundreds of lines, whose exact content we provide in the supplementary material.
We content ourselves here with an outline of their purpose.
As with prompts, and unlike in traditional programming, the relation between the text of a skill and its effect is not deterministic.
Skills are therefore benchmarked for their effectiveness, which can vary across models.
No rigorous engineering discipline for agent skills exists at present.

\paragraph{LaTeX Style Guide}

Despite its name, this skill is more than a typesetting guide.
Its rules encode how mathematical content is to be structured and articulated, so that the structure becomes legible in the TeX source and thereby visible to the model.

\paragraph{Citation Summaries}

The first paragraph of this skill reads:
\begin{prompt}
Citation summaries bridge the manuscript's `\textbackslash{}cite\{\}` claims and the original source. They are a quality-control artifact used to verify that every cited statement in the manuscript faithfully represents what the original source actually says.
\end{prompt}

For each cited source, this skill produces a summary that transfers the source's notation to that of the manuscript and states its assumptions and conclusions explicitly.
The summaries form a quality-control bridge between a citation and its source, and they put the cited statements in a form the model can handle directly.

\paragraph{Testing by Computation}

The first paragraph of this skill reads as follows:
\begin{prompt}
The manuscript may be accompanied by verification in `code` that is written in Julia.
\end{prompt}

Formulas and simplifications are checked by short computer programs, here written in \texttt{Julia}; a reader who prefers another system can adapt them, for instance to \texttt{Sage}.
In the sessions this skill was called the verification skill; we renamed it to sharpen the contrast with formal verification.

\end{appendix}

\vspace{1.5\baselineskip}
\phantomsection
\addcontentsline{toc}{section}{References}
\markright{References}
\label{sec:references}
{
  \sloppy
  \linespread{0.8}
  \printbibliography[heading=none]
}

\filbreak
\Needspace*{5\baselineskip}
\noindent%
\rule{\textwidth}{0.15em}
\\\nopagebreak

{\small\noindent
  Matthew Ortiz\\\nopagebreak
  Department of Mathematics\\\nopagebreak
  University of North Texas\\\nopagebreak
  Denton, TX 76203, USA\\\nopagebreak
  E-mail: \url{matthewortiz2@my.unt.edu}
}\vspace{.5\baselineskip}

{\small\noindent
  Martin Raum\\\nopagebreak
  Department of Mathematical Sciences\\\nopagebreak
  Chalmers University of Technology and University of Gothenburg\\\nopagebreak
  SE-412 96 Gothenburg, Sweden\\\nopagebreak
  E-mail: \url{martin@raum-brothers.eu}
}\vspace{.5\baselineskip}

{\small\noindent
  Olav K. Richter\\\nopagebreak
  Department of Mathematics\\\nopagebreak
  University of North Texas\\\nopagebreak
  Denton, TX 76203, USA\\\nopagebreak
  E-mail: \url{richter@unt.edu}%
}%

\ifdraft{%
  \listoftodos%
}

\end{document}